\documentclass[11pt]{article}
\usepackage{amsfonts}
\usepackage{latexsym}
\usepackage{cite}
\usepackage{amsmath,amsfonts,latexsym,amssymb}
\usepackage[mathscr]{eucal}
\usepackage{cases,color}
\usepackage{amsthm}

\usepackage[bf,small]{caption2}
\usepackage{float}
\usepackage{graphicx}
\usepackage{amsmath}
\usepackage{amssymb}
\usepackage[all]{xy}

\newtheorem{theorem}{theorem}[section]
\newtheorem{lemma}[theorem]{Lemma}
\newtheorem{remark}[theorem]{Remark}
\newtheorem{thm}[theorem]{Theorem}

\begin{document}

\title{\vspace{-2cm}\textbf{The ${\rm SL}(2,\mathbb{C})$-character variety of the Borromean link}}
\author{\Large Haimiao Chen \qquad  Tiantian Yu}
\date{}
\maketitle

\begin{abstract}
  For the Borromean link, we determine its irreducible ${\rm SL}(2,\mathbb{C})$-character variety, and find a formula for the twisted Alexander polynomial as a function on the character variety.

\end{abstract}

\section{Introduction}

Fix $G={\rm SL}(2,\mathbb{C})$. Let $\Gamma$ be a finitely presented group. A homomorphism $\rho:\Gamma\to G$ is called a {\it $G$-representation} of $\Gamma$. Call $\rho$ {\it reducible} if elements of ${\rm Im}(\rho)$ have a common eigenvector; otherwise, call $\rho$ {\it irreducible}. The {\it $G$-representation variety} is defined as $\mathcal{R}(\Gamma):=\hom(\Gamma,G)$. The subset $\mathcal{R}^{\rm irr}(\Gamma)$ consisting of irreducible representations form an open subspace of $\mathcal{R}(\Gamma)$.

For $\rho\in\mathcal{R}(\Gamma)$, its {\it character} is defined as the function $\chi_\rho:\Gamma\to\mathbb{C}$ sending $x$ to ${\rm tr}(\rho(x))$. Call $\mathcal{X}(\Gamma):=\{\chi_\rho\colon\rho\in\mathcal{R}(\Gamma)\}$ the $G$-{\it character variety} of $\Gamma$.
As a well-known fact, two irreducible representations $\rho,\rho'$ are conjugate (meaning that there exists $\mathbf{a}\in G$ such that $\rho'(x)=\mathbf{a}\rho(x)\mathbf{a}^{-1}$ for all $x\in\Gamma$) if and only if $\chi_\rho=\chi_{\rho'}$. Thus, an irreducible representation with conjugacy ambiguity fixed can be identified with its character. We focus on the {\it irreducible $G$-character variety} $\mathcal{X}^{\rm irr}(\Gamma)=\{\chi_\rho\colon\rho\in\mathcal{R}^{\rm irr}(\Gamma)\}.$

When $\Gamma=\pi(L):=\pi_1(S^3\setminus L)$ for a link $L$, we call $\mathcal{X}^{\rm irr}(L):=\mathcal{X}^{\rm irr}(\pi(L))$ the (irreducible) $G$-character variety of $L$.
We shall flexibly switch between characters and representations. Computing with representations are more convenient, but often needs non-canonical choices; working on characters removes the conjugacy redundance, so as to keep the description neat.

Studying the character variety of a link is a manner of understanding the link group.
In the literature, there have been various results on character varieties of knots. However, concrete computations for links are still rare. Twisted Whitehead links and $(-2,2m+1,2n)$-pretzel links were studied in \cite{Tr16}; twisted Hopf links were studied in \cite{GM23}.
See \cite{PP20} and the references in \cite{PP20,GM23} for recently new achievements about character varieties.

Although in principle the character variety for any finitely presented group can be determined algorithmically \cite{ABL18}, sometimes it is still worth finding the character variety ``by hand".
In \cite{Ch22} we proposed an easy method for investigating two elements $\mathbf{a},\mathbf{b}\in{\rm SL}(2,\mathbb{C})$ with ${\rm tr}(\mathbf{a})={\rm tr}(\mathbf{b})$ and $\mathbf{a}\mathbf{b}$ fixed. In this paper we enhance it and develop new techniques, to extract much information on $\mathbf{a},\mathbf{b}$ when the commutator $\mathbf{a}\mathbf{b}\mathbf{a}^{-1}\mathbf{b}^{-1}$ is prescribed. These are applied to determine the irreducible ${\rm SL}(2,\mathbb{C})$-character variety for the Borromean link (denoted by $B$).

The Borromean link, also called the Borromean rings, is a very interesting link. It is alternating, is the Brunnian link with the minimal number of crossings, is universal \cite{HLM83}, and was one of the earliest examples proved to be hyperbolic \cite{Ri79}.

Recently, Miura and Suzuki \cite{MS22} determined the skein algebra of $S^3\setminus B$, which modulo its nilradical is isomorphic to the coordinate ring of the whole character ${\rm SL}(2,\mathbb{C})$-variety of $\pi(B)$.
Our approach is more elementary and looks into the geometric structure of the character variety. We are able to give a clear decomposition of $\mathcal{X}^{\rm irr}(B)$ such that each part has a distinguished feature.
In particular, the {\it canonical component} admits a nice description, and is exhibited as a $2$-fold regular cover over the complement of three hypersurfaces in $(\mathbb{C}^\ast)^3$.

As another important contribution, we derive a formula for the {\it twisted Alexander polynomial} (TAP) associated to each irreducible representation; being conjugacy invariant, TAP can be regarded as a function on the irreducible character variety. The {\it twisted Reidemeister torsion}, which is also an interesting invariant, is determined by TAP in a simple way (see \cite{NT19} Theorem 2.2).

Benefiting from our results, for 3-manifolds $M$ resulting from surgeries on $B$ (see \cite{VMT01} for examples), a description of $\mathcal{X}^{\rm irr}(\pi_1(M))$ will be induced from that of $\mathcal{X}^{\rm irr}(B)$, and similarly as in \cite{NT19}, the twisted Reidemeister torsion associated to a representation $\pi_1(M)\to{\rm SL}(2,\mathbb{C})$ can be obtained by the gluing formula presented in \cite{Jo}.

\section{Algebraic techniques}

In this paper, we use $\sqrt{-1}$ (instead of $i$) to denote the imaginary unit.

For an element $x$ of a group, we often denote $\overline{x}$ for $x^{-1}$; no confusion will arise, as complex conjugate will never come into play.
Denote ${\rm Cen}(x)$ for the centralizer of $x$.
Given $x,y$, denote $x\lrcorner y$ for $xy\overline{x}$, and $[x,y]$ for the commutator $xy\overline{x}\overline{y}$. Be careful that $\overline{x}\overline{y}=x^{-1}y^{-1}$, not $(xy)^{-1}$.

Let $\mathcal{D},\mathcal{T}_+,\mathcal{T}_-$ respectively denote the subgroup of ${\rm SL}(2,\mathbb{C})$ consisting of diagonal, upper-triangular, lower-triangular matrices.

We use boldface letters to denote elements of ${\rm SL}(2,\mathbb{C})$. Let
\begin{align*}
\mathbf{w}=\left(\begin{array}{cc} 0 & 1 \\ -1 & 0 \end{array}\right),  \qquad
\mathbf{d}(\kappa)=\left(\begin{array}{cc} \kappa & 0 \\ 0 & \kappa^{-1} \end{array}\right), \qquad
\mathbf{p}(u)=\left(\begin{array}{cc} 1 & u \\ 0 & 1 \end{array}\right).
\end{align*}
Let $\mathbf{e}$ denote the identity matrix, i.e. $\mathbf{e}=\mathbf{d}(1)=\mathbf{p}(0)$; let $\mathbf{p}=\mathbf{p}(1)$.

For $\mathbf{x}\in{\rm SL}(2,\mathbb{C})$, let $\mathbf{x}_{ij}$ denote its $(i,j)$-entry.

Remember that, by Hamilton-Cayley Theorem,
\begin{align}
\mathbf{x}^2={\rm tr}(\mathbf{x})\cdot\mathbf{x}-\mathbf{e}.  \label{eq:H-C}
\end{align}

For $\lambda+\lambda^{-1}\ne t^2-2,\pm 2$ and $\mu\ne 0$, put
\begin{align*}
\mathbf{h}_{t}^{\lambda}(\mu)=
\frac{1}{\lambda+1}\left(\begin{array}{cc} \lambda t & \mu \\ \delta\lambda\mu^{-1} & t \end{array}\right), \qquad \delta=t^2-\lambda-\lambda^{-1}-2.
\end{align*}
For $t\ne 0$, put
\begin{align*}
\mathbf{k}_{t}(\alpha)=\left(\begin{array}{cc} \alpha+t/2  & (2t)^{-1}(t^2/4-1-\alpha^2) \\ 2t & -\alpha+t/2  \end{array}\right).
\end{align*}
The reason for introducing these matrices is clear from the following lemmas. 

\begin{lemma}\label{lem:key}
Suppose $\mathbf{a},\mathbf{b}\in{\rm SL}(2,\mathbb{C})$ with ${\rm tr}(\mathbf{a})={\rm tr}(\mathbf{b})=t$.
\begin{enumerate}
  \item[\rm(a)] If $\mathbf{a}\mathbf{b}=\mathbf{d}(\lambda)$ with $\lambda+\lambda^{-1}\ne t^2-2, \pm 2$,
                then there exists $\mu\ne 0$ such that $\mathbf{a}=\mathbf{h}_{t}^{\lambda}(\mu)$ and $\mathbf{b}=\mathbf{h}_{t}^{\lambda}(-\lambda^{-1}\mu)$.
  \item[\rm(a')] If $\mathbf{a}\mathbf{b}=\mathbf{d}(\lambda)$ with $\lambda+\lambda^{-1}=t^2-2\ne\pm 2$,
                then $\mathbf{a}, \mathbf{b}\in \mathcal{T}_+$ or $\mathbf{a}, \mathbf{b}\in \mathcal{T}_-$.
  \item[\rm(b)] If $\mathbf{a}\mathbf{b}=\mathbf{p}$, then $\mathbf{a}, \mathbf{b}\in \mathcal{T}_+$.
  \item[\rm(c)] If $\mathbf{a}\mathbf{b}=-\mathbf{p}$ and $t\ne 0$, then $\mathbf{a}=\mathbf{k}_{t}(\alpha)$ and $\mathbf{b}=\mathbf{k}_{t}(\alpha-t)$ for some $\alpha$.
  \item[\rm(c')] If $\mathbf{a}\mathbf{b}=-\mathbf{p}$ and $t=0$, then $\mathbf{a}, \mathbf{b}\in \mathcal{T}_+$.
\end{enumerate}
\end{lemma}
\begin{proof}
(a),(b),(c) were given in \cite{Ch22} Lemma 2.2. We prove (a), (a'), (c') here, and refer to \cite{Ch22} for (b), (c).

To show (a) and (a'), assume $\lambda\ne\pm1$, and suppose $\mathbf{a}=(a_{ij})_{2\times 2}$, $\mathbf{b}=(b_{ij})_{2\times 2}$, with
$a_{11}a_{22}-a_{12}a_{21}=b_{11}b_{22}-b_{12}b_{21}=1.$
Then $\mathbf{a}\mathbf{b}=\mathbf{d}(\lambda)$ implies $\mathbf{b}=\overline{\mathbf{a}}\mathbf{d}(\lambda)$, so that $b_{11}=\lambda a_{22}$, $b_{22}=\lambda^{-1}a_{11}$. Since $t=a_{11}+a_{22}$ and $t=b_{11}+b_{22}=\lambda a_{22}+\lambda^{-1}a_{11}$, we have
$a_{11}=\lambda t/(\lambda+1)$, $a_{22}=t/(\lambda+1)$.
Now if $\lambda+\lambda^{-1}\ne t^2-2$, then
$$a_{12}a_{21}=a_{11}a_{22}-1=\frac{\lambda(t^2-2-\lambda-\lambda^{-1})}{(\lambda+1)^2}\ne 0,$$
hence $\mathbf{a}=\mathbf{h}_{t}^{\lambda}(a_{12})$; it follows that $\mathbf{b}=\mathbf{h}_{t}^{\lambda}(-\lambda^{-1}a_{12})$.
If $\lambda+\lambda^{-1}=t^2-2$, then $a_{12}a_{21}=0$, implying $\mathbf{a}\in\mathcal{T}_+$ or $\mathbf{a}\in\mathcal{T}_-$; respectively, $\mathbf{b}\in\mathcal{T}_+$ or $\mathbf{b}\in\mathcal{T}_-$.

For (c'), just note $(-\mathbf{a})\mathbf{b}=\mathbf{p}$ and ${\rm tr}(-\mathbf{a})={\rm tr}(\mathbf{b})=0$, and apply (b).
\end{proof}

\begin{lemma}\label{lem:commutator}
Suppose ${\rm tr}(\mathbf{a}_i)=t_i$, $i=1,2$. Let $\mathbf{g}=[\mathbf{a}_1,\mathbf{a}_2]$.
\begin{enumerate}
  \item[\rm(a)] If $\lambda\ne\pm 1$ and $\delta_i:=t_i^2-\lambda-\lambda^{-1}-2\ne 0$, $i=1,2$, then $\mathbf{g}=\mathbf{d}(\lambda)$ if and
        only if $\mathbf{a}_1=\mathbf{h}_{t_1}^\lambda(\mu)$, $\overline{\mathbf{a}_2}=\mathbf{h}_{t_2}^\lambda(\nu)$ for some $\mu,\nu\ne 0$ with
        \begin{align}
        (\lambda-1)t_1t_2=\lambda\delta_1\mu^{-1}\nu-\delta_2\mu\nu^{-1}. \label{eq:condition-1}
        \end{align}
  \item[\rm(b)] If $\mathbf{g}=\mathbf{p}$, then $\mathbf{a}_1,\mathbf{a}_2\in\mathcal{T}_+$.
  \item[\rm(c)] If $t_1t_2\ne 0$, then $\mathbf{g}=-\mathbf{p}$ if and only if $\mathbf{a}_1=\mathbf{k}_{t_1}(\alpha)$, $\overline{\mathbf{a}_2}=\mathbf{k}_{t_2}(\beta)$ for some $\alpha,\beta$ with
        \begin{align}
        t_1^{-1}t_2(\alpha^2+1)+t_1t_2^{-1}(\beta^2+1)=2\alpha\beta+t_2\alpha-t_1\beta.    \label{eq:condition-2}
        \end{align}
  \item[\rm(d)] $\mathbf{g}=-\mathbf{e}$ if and only if $t_1=t_2={\rm tr}(\mathbf{a}_1\mathbf{a}_2)=0$.
  \item[\rm(e)] ${\rm tr}(\mathbf{g})=2$ if and only if $\mathbf{a}_1,\mathbf{a}_2$ have a common eigenvector.
\end{enumerate}
\end{lemma}

\begin{proof}
(a) ($\Rightarrow$) Suppose $\mathbf{g}=\mathbf{d}(\lambda)$. Since $\mathbf{a}_1\cdot\mathbf{a}_2\lrcorner\overline{\mathbf{a}_1}=\mathbf{g}=\mathbf{d}(\lambda)$ and
${\rm tr}(\mathbf{a}_1)={\rm tr}(\mathbf{a}_2\lrcorner\overline{\mathbf{a}_1})=t_1$, by Lemma \ref{lem:key} (a), there exists $\mu\ne 0$ such that $\mathbf{a}_1=\mathbf{h}_{t_1}^{\lambda}(\mu)$ and $\mathbf{a}_2\lrcorner\overline{\mathbf{a}_1}=\mathbf{h}_{t_1}^{\lambda}(-\lambda^{-1}\mu)$.
Since $\mathbf{a}_1\lrcorner\mathbf{a}_2\cdot\overline{\mathbf{a}_2}=\mathbf{d}(\lambda)$ and ${\rm tr}(\mathbf{a}_1\lrcorner\mathbf{a}_2)={\rm tr}(\overline{\mathbf{a}_2})=t_2$, by Lemma \ref{lem:key} (a) again,
there exists $\nu\ne 0$ such that $\mathbf{a}_1\lrcorner\mathbf{a}_2=\mathbf{h}_{t_2}^{\lambda}(-\lambda\nu)$ and $\overline{\mathbf{a}_2}=\mathbf{h}_{t_2}^{\lambda}(\nu)$.
Note that $\overline{\mathbf{a}_1}\overline{\mathbf{a}_2}=\overline{\mathbf{a}_2}\cdot\mathbf{a}_2\lrcorner\overline{\mathbf{a}_1}$ requires
\begin{align}
\mathbf{h}_{t_1}^{\lambda}(\mu)^{-1}\mathbf{h}_{t_2}^{\lambda}(\nu)=\mathbf{h}_{t_2}^{\lambda}(\nu)\mathbf{h}_{t_1}^{\lambda}(-\lambda^{-1}\mu).
\label{eq:consistence-1}
\end{align}
By direct computation,
\begin{align*}
\mathbf{h}_{t_1}^{\lambda}(\mu)^{-1}\mathbf{h}_{t_2}^{\lambda}(\nu)
&=\frac{1}{(\lambda+1)^2}\left(\begin{array}{cc} \lambda(t_1t_2-\delta_2\mu\nu^{-1}) & t_1\nu-t_2\mu \\
\lambda^2(\delta_2t_1\nu^{-1}-\delta_1t_2\mu^{-1}) & \lambda(t_1t_2-\delta_1\nu\mu^{-1}) \end{array}\right),   \\
\mathbf{h}_{t_2}^{\lambda}(\nu)\mathbf{h}_{t_1}^{\lambda}(-\lambda^{-1}\mu)
&=\frac{1}{(\lambda+1)^2}\left(\begin{array}{cc} \lambda^2(t_1t_2-\delta_1\nu\mu^{-1}) & t_1\nu-t_2\mu \\
\lambda^2(\delta_2t_1\nu^{-1}-\delta_1t_2\mu^{-1}) & t_1t_2-\delta_2\mu\nu^{-1} \end{array}\right).
\end{align*}
Hence (\ref{eq:consistence-1}) is equivalent to (\ref{eq:condition-1}).

($\Leftarrow$) If $\mathbf{a}_1=\mathbf{h}_{t_1}^\lambda(\mu)$, $\overline{\mathbf{a}_2}=\mathbf{h}_{t_2}^\lambda(\nu)$ for some $\mu,\nu\ne 0$ satisfying (\ref{eq:condition-1}), then
\begin{align*}
\mathbf{g}&=\mathbf{h}_{t_1}^\lambda(\mu)\mathbf{h}_{t_2}^\lambda(\nu)^{-1}\cdot\mathbf{h}_{t_1}^\lambda(\mu)^{-1}\mathbf{h}_{t_2}^\lambda(\nu)
=\mathbf{h}_{t_1}^\lambda(\mu)\mathbf{h}_{t_2}^\lambda(\nu)^{-1}\cdot\mathbf{h}_{t_2}^{\lambda}(\nu)\mathbf{h}_{t_1}^{\lambda}(-\lambda^{-1}\mu)  \\
&=\mathbf{h}_{t_1}^\lambda(\mu)\mathbf{h}_{t_1}^{\lambda}(-\lambda^{-1}\mu)=\mathbf{d}(\lambda).
\end{align*}

(b) Suppose $\mathbf{g}=\mathbf{p}$. Since $\mathbf{a}_1\cdot\mathbf{a}_2\lrcorner\overline{\mathbf{a}_1}=\mathbf{p}$ and ${\rm tr}(\mathbf{a}_1)={\rm tr}(\mathbf{a}_2\lrcorner\overline{\mathbf{a}_1})$, by Lemma \ref{lem:key} (b), $\mathbf{a}_1, \mathbf{a}_2\lrcorner\overline{\mathbf{a}_1}\in\mathcal{T}_+$.
Since $\mathbf{a}_1\lrcorner\mathbf{a}_2\cdot\overline{\mathbf{a}_2}=\mathbf{p}$ and ${\rm tr}(\mathbf{a}_1\lrcorner\mathbf{a}_2)={\rm tr}(\overline{\mathbf{a}_2})$, by Lemma \ref{lem:key} (b) again, $\mathbf{a}_1\lrcorner\mathbf{a}_2,\overline{\mathbf{a}_2}\in\mathcal{T}_+$. Hence $\mathbf{a}_1,\mathbf{a}_2\in\mathcal{T}_+$.

(c) ($\Rightarrow$) Suppose $\mathbf{g}=-\mathbf{p}$. Since $\mathbf{a}_1\cdot\mathbf{a}_2\lrcorner\overline{\mathbf{a}_1}=-\mathbf{p}$ and ${\rm tr}(\mathbf{a}_1)={\rm tr}(\mathbf{a}_2\lrcorner\overline{\mathbf{a}_1})=t_1$, by Lemma \ref{lem:key} (c), there exists $\alpha$ such that $\mathbf{a}_1=\mathbf{k}_{t_1}(\alpha)$ and $\mathbf{a}_2\lrcorner\overline{\mathbf{a}_1}=\mathbf{k}_{t_1}(\alpha-t_1)$.
Since $\mathbf{a}_1\lrcorner\mathbf{a}_2\cdot\overline{\mathbf{a}_2}=-\mathbf{p}$ and ${\rm tr}(\mathbf{a}_1\lrcorner\mathbf{a}_2)={\rm tr}(\overline{\mathbf{a}_2})=t_2$, by Lemma \ref{lem:key} (c) again,
there exists $\beta$ such that $\mathbf{a}_1\lrcorner\mathbf{a}_2=\mathbf{k}_{t_2}(\beta+t_2)$ and $\overline{\mathbf{a}_2}=\mathbf{k}_{t_2}(\beta)$.
Note that $\overline{\mathbf{a}_1}\overline{\mathbf{a}_2}=\overline{\mathbf{a}_2}\cdot\mathbf{a}_2\lrcorner\overline{\mathbf{a}_1}$ requires
\begin{align}
\mathbf{k}_{t_1}(\alpha)^{-1}\mathbf{k}_{t_2}(\beta)=\mathbf{k}_{t_2}(\beta)\mathbf{k}_{t_1}(\alpha-t_1).   \label{eq:consistence-2}
\end{align}
By direct computation, $\mathbf{k}_{t_1}(\alpha)^{-1}\mathbf{k}_{t_2}(\beta)$ and $\mathbf{k}_{t_2}(\beta)\mathbf{k}_{t_1}(\alpha-t_1)$ equals
\begin{align*}
&\left(\begin{array}{cc} \frac{1}{2}(t_1\beta-t_2\alpha)-\alpha\beta+t_1^{-1}t_2(\alpha^2+1) & \ast \\
2(t_2\alpha-t_1\beta) & \frac{1}{2}(t_2\alpha-t_1\beta)-\alpha\beta+t_1t_2^{-1}(\beta^2+1) \end{array}\right),   \\
&\left(\begin{array}{cc} \frac{1}{2}(t_2\alpha-t_1\beta)+\alpha\beta-t_1t_2^{-1}(\beta^2+1) & \ast' \\
2(t_2\alpha-t_1\beta) & \frac{3}{2}(t_2\alpha-t_1\beta)+\alpha\beta-t_1^{-1}t_2(\alpha^2+1) \end{array}\right),
\end{align*}
respectively, where $\ast,\ast'$ are the numbers determined by that the matrices have determinant $1$, and are not irrelevant to us.
So (\ref{eq:consistence-2}) is equivalent to (\ref{eq:condition-2}).

($\Leftarrow$) If $\mathbf{a}_1=\mathbf{k}_{t_1}(\alpha)$, $\overline{\mathbf{a}_2}=\mathbf{k}_{t_2}(\beta)$ for some $\alpha,\beta$ satisfying (\ref{eq:condition-2}), then
\begin{align*}
\mathbf{g}&=\mathbf{k}_{t_1}(\alpha)\mathbf{k}_{t_2}(\beta)^{-1}\cdot\mathbf{k}_{t_1}(\alpha)^{-1}\mathbf{k}_{t_2}(\beta)
=\mathbf{k}_{t_1}(\alpha)\mathbf{k}_{t_2}(\beta)^{-1}\cdot\mathbf{k}_{t_2}(\beta)\mathbf{k}_{t_1}(\alpha-t_1)  \\
&=\mathbf{k}_{t_1}(\alpha)\mathbf{k}_{t_1}(\alpha-t_1)=-\mathbf{p}.
\end{align*}

(d) For each $\mathbf{x}\in{\rm SL}(2,\mathbb{C})$, by (\ref{eq:H-C}) we have ${\rm tr}(\mathbf{x})=0\Leftrightarrow\mathbf{x}^2=-\mathbf{e}$.

($\Leftarrow$) If $t_1=t_2={\rm tr}(\mathbf{a}_1\mathbf{a}_2)=0$, then $\mathbf{g}=\mathbf{a}_1\mathbf{a}_2(-\mathbf{a}_1)(-\mathbf{a}_2)=(\mathbf{a}_1\mathbf{a}_2)^2=-\mathbf{e}$.

($\Rightarrow$) Suppose $\mathbf{g}=-\mathbf{e}$. Then $-\mathbf{a}_1=\mathbf{a}_2\lrcorner\mathbf{a}_1$ and $-\mathbf{a}_2=\mathbf{a}_1\lrcorner\mathbf{a}_2$, respectively implying $-t_1=t_1$ and $-t_2=t_2$. So $t_1=t_2=0$. Now $\mathbf{a}_1^2=\mathbf{a}_2^2=-\mathbf{e}$, hence $-\mathbf{e}=\mathbf{g}=\mathbf{a}_1\mathbf{a}_2(-\mathbf{a}_1)(-\mathbf{a}_2)=(\mathbf{a}_1\mathbf{a}_2)^2$, implying ${\rm tr}(\mathbf{a}_1\mathbf{a}_2)=0$.

(e) This is a part of \cite{Go09} Proposition 2.3.1. We reprove it for completeness.

($\Leftarrow$) If $\mathbf{a}_1,\mathbf{a}_2$ have a common eigenvector, then there exists $\mathbf{c}$ such that $\mathbf{c}\lrcorner\mathbf{a}_1,\mathbf{c}\lrcorner\mathbf{a}_2\in\mathcal{T}_+$, hence
${\rm tr}(\mathbf{g})={\rm tr}([\mathbf{c}\lrcorner\mathbf{a}_1,\mathbf{c}\lrcorner\mathbf{a}_2])=2$.

($\Rightarrow$) Suppose ${\rm tr}(\mathbf{g})=2$. If $\mathbf{g}=\mathbf{e}$, then $\mathbf{a}_1\mathbf{a}_2=\mathbf{a}_2\mathbf{a}_1$, so $\mathbf{a}_1,\mathbf{a}_2$ have a common eigenvector. If $\mathbf{g}\ne\mathbf{e}$, then up to conjugacy we may assume $\mathbf{g}=\mathbf{p}$; by (b), $\mathbf{a}_1,\mathbf{a}_2\in\mathcal{T}_+$, implying that $\mathbf{a}_1,\mathbf{a}_2$ have a common eigenvector.
\end{proof}

\section{The character variety}

\begin{figure}[h]
  \centering
  \includegraphics[width=8cm]{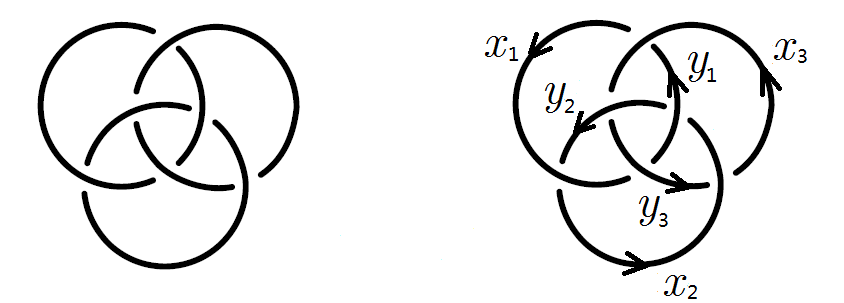} \\
  \caption{Left: the Borromean link $B$. Right: a direction is chosen for each arc.} \label{fig:Borromean}
\end{figure}

A symmetric diagram for the Borromean link $B$ is given in Figure \ref{fig:Borromean}.

The Wirtinger presentation gives
\begin{align}
\pi(B)\cong\langle x_1,x_2,x_3,y_1,y_2,y_3\mid\overline{y_1}x_3\lrcorner x_1, \overline{y_2}x_1\lrcorner x_2, \overline{y_3}x_2\lrcorner x_3, \overline{x_3}\overline{y_2}\lrcorner y_3, \overline{x_2}\overline{y_1}\lrcorner y_2\rangle;  \label{eq:Wirtinger}
\end{align}
with $y_i$'s substituted, this can be transformed into
\begin{align}
\pi(B)\cong\langle x_1,x_2,x_3\mid [x_2,[x_3,\overline{x_1}]],\ [x_3,[x_1,\overline{x_2}]]\rangle.   \label{eq:presentation}
\end{align}
As is always true in Wirtinger presentation, in (\ref{eq:Wirtinger}), the ``omitted relation" $\overline{x_1}\overline{y_3}\lrcorner y_1$ follows from the five given relations.
So $[x_2,[x_3,\overline{x_1}]]$ and $[x_3,[x_1,\overline{x_2}]]$ actually imply $[x_1,[x_2,\overline{x_3}]]$.
To be explicit, assuming $x_2[x_3,\overline{x_1}]=[x_3,\overline{x_1}]x_2$ and $x_3[x_1,\overline{x_2}]=[x_1,\overline{x_2}]x_3$, one can deduce
$x_1[x_2,\overline{x_3}]=[x_2,\overline{x_3}]x_1$ as follows: $x_3[x_1,\overline{x_2}]=[x_1,\overline{x_2}]x_3$ implies $x_2\overline{x_3}\overline{x_2}=x_1x_2\overline{x_1}\overline{x_3}x_1\overline{x_2}\overline{x_1}$, so
\begin{align*}
[x_2,\overline{x_3}]&=x_2\overline{x_3}\overline{x_2}\cdot x_3=x_1x_2\overline{x_1}\overline{x_3}x_1\overline{x_2}\overline{x_1}\cdot x_3
=x_1x_2\overline{x_3}\cdot[x_3,\overline{x_1}]\overline{x_2}\cdot\overline{x_1}x_3 \\
&=x_1x_2\overline{x_3}\cdot\overline{x_2}[x_3,\overline{x_1}]\cdot\overline{x_1}x_3=x_1[x_2,\overline{x_3}]\overline{x_1}.
\end{align*}

The element represented by the longitude paired with $x_i$ is $[x_{i+1},\overline{x_{i-1}}]$.
The subscripts are understood as modulo 3, as done throughout this paper.

Given $\mathbf{x}_1,\mathbf{x}_2,\mathbf{x}_3\in{\rm SL}(2,\mathbb{C})$, there exists a (unique) representation
$\rho:\pi(B)\to{\rm SL}(2,\mathbb{C})$ with $\rho(x_i)=\mathbf{x}_i$ if and only if $[\mathbf{x}_i,[\mathbf{x}_{i+1},\overline{\mathbf{x}_{i-1}}]]=\mathbf{e}$ for $i=1,2,3$.
As explained above, each of the three equations follows from the other two.
We always assume that $\rho$ is irreducible; equivalently, $\mathbf{x}_1,\mathbf{x}_2,\mathbf{x}_3$ have no common eigenvector.

If $\mathbf{x}_i\in\{\pm\mathbf{e}\}$, then $\mathbf{x}_{i-1},\mathbf{x}_{i+1}$ can be arbitrary.

From now on, we assume $\mathbf{x}_i\ne\pm\mathbf{e}$ for all $i$.

Let $\mathbf{g}_i=[\mathbf{x}_{i+1},\overline{\mathbf{x}_{i-1}}]$.
Let $\lambda_i+\lambda_i^{-1}={\rm tr}(\mathbf{g}_i)$, $t_i={\rm tr}(\mathbf{x}_i)$; let $t_{ij}=t_{i,j}={\rm tr}(\mathbf{x}_i\mathbf{x}_j)$ for $i\ne j$; let $t_{123}={\rm tr}(\mathbf{x}_1\mathbf{x}_2\mathbf{x}_3)$.

We shall frequently use the following facts which can be directly verified:
\begin{alignat*}{2}
{\rm Cen}(\mathbf{d}(\kappa))&=\mathcal{D}, \qquad \kappa\ne\pm1;  \\
{\rm Cen}(\mathbf{p})={\rm Cen}(-\mathbf{p})&=\{\epsilon\mathbf{p}(u)\colon \epsilon\in\{\pm1\},\ u\in\mathbb{C}\}.
\end{alignat*}
Stated alternatively, if $\mathbf{a}=\mathbf{d}(\kappa)$ or $\mathbf{a}\in\{\pm\mathbf{p}\}$, then elements of ${\rm Cen}(\mathbf{a})$ are of the form $\mu\mathbf{a}+\nu\mathbf{e}$; since this property is invariant under conjugation, it holds for all $\mathbf{a}\ne\pm\mathbf{e}$.
As a consequence, ${\rm Cen}(\mathbf{a})\subset\mathcal{T}_+$ for any $\mathbf{a}\in\mathcal{T}_+\setminus\{\pm\mathbf{e}\}$.

\subsection{$\mathbf{g}_i=\pm\mathbf{e}$ for some $i$}

\begin{lemma} \label{lem:g=e}
{\rm(i)} $\mathbf{g}_i=\mathbf{e}$ if and only if $\lambda_i=1$.

{\rm(ii)} Suppose $\mathbf{g}_i=\mathbf{e}$. Then up to conjugacy $\mathbf{x}_i=\mathbf{w}$, $\mathbf{x}_{i-1}=\mathbf{d}(\kappa_{i-1})$, $\mathbf{x}_{i+1}=\mathbf{d}(\kappa_{i+1})$, for some $\kappa_{i\pm1}$ with $\kappa_{i\pm1}+\kappa_{i\pm1}^{-1}=t_{i\pm1}\ne\pm 2$.
\end{lemma}

\begin{proof}
(i) Clearly, $\mathbf{g}_i=\mathbf{e}$ implies $\lambda_i=1$.

If $\lambda_i=1$ but $[\mathbf{x}_{i+1},\overline{\mathbf{x}_{i-1}}]=\mathbf{g}_i\ne\mathbf{e}$, then up to conjugacy we may assume $\mathbf{g}_i=\mathbf{p}$, so that $\mathbf{x}_i\in{\rm Cen}(\mathbf{p})\subset\mathcal{T}_+$. By Lemma \ref{lem:commutator} (b), $\mathbf{x}_{i-1},\mathbf{x}_{i+1}\in\mathcal{T}_+$, contradicting the irreducibility. Hence $\lambda_i=1$ implies $\mathbf{g}_i=\mathbf{e}$,

(ii) We have $\mathbf{x}_{i+1}\in{\rm Cen}(\mathbf{x}_{i-1})$, and also
$\mathbf{x}_i\lrcorner\mathbf{x}_{i+1}=\mathbf{x}_{i+1}\overline{\mathbf{g}_{i-1}}\in{\rm Cen}(\mathbf{x}_{i-1})$.

If $t_{i-1}=2\epsilon$ with $\epsilon\in\{\pm 1\}$, then up to conjugacy we may assume $\mathbf{x}_{i-1}=\epsilon\mathbf{p}$, so that $\mathbf{x}_{i+1}=\epsilon'\mathbf{p}(v)$ for some $\epsilon'\in\{\pm1\}$, $v\ne 0$.
However, one can verify that $(\mathbf{x}_i\lrcorner\mathbf{x}_{i+1})_{21}=-\epsilon'v((\mathbf{x}_i)_{21})^2$, so
$\mathbf{x}_i\lrcorner\mathbf{x}_{i+1}\in{\rm Cen}(\mathbf{x}_{i-1})$ would force $\mathbf{x}_i\in\mathcal{T}_+$.
This contradicts the irreducibility.

Thus, $t_{i-1}\ne\pm 2$. Up to conjugacy we can assume $\mathbf{x}_{i-1}=\mathbf{d}(\kappa_{i-1})$ with $\kappa_{i-1}+\kappa_{i-1}^{-1}=t_{i-1}$, so that $\mathbf{x}_{i+1}\in\mathcal{D},\mathbf{x}_i\lrcorner\mathbf{x}_{i+1}\in\mathcal{D}$.
Due to the assumption $\mathbf{x}_{i+1}\ne\pm\mathbf{e}$, we have $\mathbf{x}_{i+1}=\mathbf{d}(\kappa_{i+1})$ for some $\kappa_{i+1}$ with $\kappa_{i+1}+\kappa_{i+1}^{-1}=t_{i+1}\ne\pm 2$.
Note that
${\rm tr}(\mathbf{x}_i\lrcorner\mathbf{x}_{i+1})=t_{i+1}$ and $\mathbf{x}_i\lrcorner\mathbf{x}_{i+1}=\mathbf{x}_{i+1}$ is impossible (otherwise $\mathbf{x}_i\in{\rm Cen}(\mathbf{x}_{i+1})=\mathcal{D}$, contradicting the irreducibility), so $\mathbf{x}_i\lrcorner\mathbf{x}_{i+1}$ must equal $\overline{\mathbf{x}_{i+1}}$, i.e. $\mathbf{x}_{i+1}\mathbf{x}_i\mathbf{x}_{i+1}=\mathbf{x}_i$.
Consequently, $(\mathbf{x}_i)_{11}=(\mathbf{x}_i)_{22}=0$. Conjugating via a suitable element of $\mathcal{D}$, we may further set $\mathbf{x}_i=\mathbf{w}$.
\end{proof}

\begin{lemma} \label{lem:g=-e}
If $\mathbf{g}_i=-\mathbf{e}$, then up to conjugacy $\mathbf{x}_{i-1}=\mathbf{w}$, $\mathbf{x}_{i+1}=\mathbf{d}(\sqrt{-1})$, $\mathbf{x}_i=(a_{jk})_{2\times 2}$ with $a_{11}a_{22}-a_{12}a_{21}=1$ and $a_{11}a_{12}+a_{21}a_{22}=0$.
\end{lemma}

\begin{proof}
Since $[\mathbf{x}_{i+1},\overline{\mathbf{x}_{i-1}}]=\mathbf{g}_i=-\mathbf{e}$, we have $[\mathbf{x}_{i-1},\mathbf{x}_{i+1}]=-\mathbf{e}$.
By Lemma \ref{lem:commutator} (d), $t_{i+1}=t_{i-1}=t_{i-1,i+1}=0$. Up to conjugacy we may assume $\mathbf{x}_{i+1}=\mathbf{d}(\sqrt{-1})$, and then deduce from $t_{i-1}=t_{i-1,i+1}=0$ that $(\mathbf{x}_{i-1})_{11}=(\mathbf{x}_{i-1})_{22}=0$. Conjugating by a suitable element of $\mathcal{D}$, we may further assume $\mathbf{x}_{i-1}=\mathbf{w}$.

Suppose $\mathbf{x}_i=(a_{jk})_{2\times 2}$ with $a_{11}a_{22}-a_{12}a_{21}=1$.
Then
$$[\mathbf{x}_{i-1},\overline{\mathbf{x}_i}]=\left(\begin{array}{cc} a_{11}^2+a_{21}^2 & a_{11}a_{12}+a_{21}a_{22} \\
a_{11}a_{12}+a_{21}a_{22} & a_{12}^2+a_{22}^2 \end{array}\right).$$
Hence $[\mathbf{x}_{i-1},\overline{\mathbf{x}_i}]\in{\rm Cen}(\mathbf{x}_{i+1})$ is equivalent to $a_{11}a_{12}+a_{21}a_{22}=0$.
\end{proof}

In the setting of the proof,
\begin{alignat*}{2}
t_i&=a_{11}+a_{22},   &t_{i,i+1}&={\rm tr}\big(\mathbf{x}_i\mathbf{d}(\sqrt{-1})\big)=(a_{11}-a_{22})\sqrt{-1}, \\
t_{i,i-1}&={\rm tr}(\mathbf{w}\mathbf{x}_i)=a_{21}-a_{12}, \ \ \
&t_{123}&={\rm tr}\big(\mathbf{w}\mathbf{x}_i\mathbf{d}(\sqrt{-1})\big)=(a_{12}+a_{21})\sqrt{-1}.
\end{alignat*}
Thus,
\begin{align*}
a_{11}&=\frac{1}{2}(t_i-t_{i,i+1}\sqrt{-1}), &&a_{12}=-\frac{1}{2}(t_{i,i-1}+t_{123}\sqrt{-1}), \\
a_{21}&=\frac{1}{2}(t_{i,i-1}-t_{123}\sqrt{-1}), &&a_{22}=\frac{1}{2}(t_i+t_{i,i+1}\sqrt{-1}),
\end{align*}
and the conclusion of the lemma can be restated as
\begin{align*}
t_{i-1}=t_{i+1}=t_{i-1,i+1}=0, \\
t_it_{123}=t_{i,i-1}t_{i,i+1}, \\
t_i^2+t_{i,i-1}^2+t_{i,i+1}^2+t_{123}^2=4.
\end{align*}

\subsection{$\mathbf{g}_i\ne\pm\mathbf{e}$ for all $i$}

By Lemma \ref{lem:g=e} (i), $\lambda_i\ne 1$ for all $i$.

\subsubsection{$\lambda_3\ne-1$} \label{sec:lambda3-general}

Let $\lambda=\lambda_3\ne\pm1$.
Up to conjugacy, we may assume $\mathbf{g}_3=\mathbf{d}(\lambda)$ and $\mathbf{x}_3=\mathbf{d}(\kappa)$ with $\kappa\ne\pm1$ and $\kappa+\kappa^{-1}=t_3$.

Let $\delta_i=t_i^2-\lambda-\lambda^{-1}-2$, $i=1,2$. 

\begin{lemma}
$\delta_1\delta_2\ne 0$.
\end{lemma}
\begin{proof}
Assume $\delta_1=0$, i.e. $\lambda+\lambda^{-1}=t_1^2-2$.
Since $\mathbf{x}_1\cdot\overline{\mathbf{x}_2}\lrcorner\overline{\mathbf{x}_1}=\mathbf{g}_3=\mathbf{d}(\lambda)$ and ${\rm tr}(\mathbf{x}_1)={\rm tr}(\overline{\mathbf{x}_2}\lrcorner\overline{\mathbf{x}_1})=t_1$, by Lemma \ref{lem:key} (a'), $\mathbf{x}_1\in\mathcal{T}_+$ or $\mathbf{x}_1\in\mathcal{T}_-$; say $\mathbf{x}_1\in\mathcal{T}_+$.
Then $\mathbf{g}_2=[\mathbf{x}_3,\overline{\mathbf{x}_1}]\in\mathcal{T}_+$.
Since $\mathbf{g}_2\ne\mathbf{e}$, we have $\mathbf{x}_2\in{\rm Cen}(\mathbf{g}_2)\subset\mathcal{T}_+$. But this contradicts the irreducibility. The situation is similar when $\mathbf{x}_1\in\mathcal{T}_-$. Thus, $\delta_1\ne 0$.

Similarly, $\delta_2\ne 0$.
\end{proof}

An application of Lemma \ref{lem:commutator} (a) to $\mathbf{a}_1=\mathbf{x}_1$ and $\mathbf{a}_2=\overline{\mathbf{x}_2}$ shows $\mathbf{x}_1=\mathbf{h}_{t_1}^{\lambda}(\mu)$, $\mathbf{x}_2=\mathbf{h}_{t_2}^{\lambda}(\nu)$ for some $\mu,\nu\ne 0$; by conjugating via $\mathbf{d}(1/\sqrt{\mu})$, we may just assume $\mathbf{x}_1=\mathbf{h}_{t_1}^{\lambda}(1)$. Hence the condition in Lemma \ref{lem:commutator} (a) reads
\begin{align}
(\lambda-1)t_1t_2=\lambda\delta_1\nu-\delta_2\nu^{-1}.  \label{eq:case1-1}
\end{align}
Applying Lemma \ref{lem:commutator} (a) to $\mathbf{a}_1=\mathbf{x}_3\lrcorner\mathbf{x}_1=\mathbf{h}_{t_1}(\kappa^2\mu)$ and $\mathbf{a}_2=\overline{\mathbf{x}_2}$, we obtain
\begin{align}
(\lambda-1)t_1t_2=\lambda\delta_1\kappa^{-2}\nu-\delta_2\kappa^2\nu^{-1}.  \label{eq:case1-2}
\end{align}

\begin{lemma}
$t_1t_2t_3\ne 0$.
\end{lemma}
\begin{proof}
If $t_1=0$, then (\ref{eq:case1-1}), (\ref{eq:case1-2}) would respectively imply
$\lambda_3\delta_1\nu-\delta_2\nu^{-1}=0$ and $\lambda\delta_1\kappa^{-2}\nu-\delta_2\kappa^2\nu^{-1}=0$; since $\lambda\delta_1\nu\ne 0$, $\delta_2\nu^{-1}\ne0$ and $\kappa^2\ne 1$, we must have $\kappa^2=-1$, so that $t_3=0$.
Moreover, $t_{13}={\rm tr}\big(\mathbf{h}_{t_1}^{\lambda}(1)\mathbf{d}(\kappa)\big)=0$. By Lemma \ref{lem:commutator} (d), $[\mathbf{x}_1,\mathbf{x}_3]=-\mathbf{e}$, which contradicts the assumption $\mathbf{g}_2\ne-\mathbf{e}$.

Thus, $t_1\ne 0$. Similarly, $t_2\ne 0$.

By (\ref{eq:case1-1}), (\ref{eq:case1-2}), $2(\lambda-1)t_1t_2=t_3(\kappa^{-1}\lambda\delta_1\nu-\kappa\delta_2\nu^{-1})$. Hence $t_3\ne 0$.
\end{proof}

As a consequence of (\ref{eq:case1-1}) and (\ref{eq:case1-2}),
\begin{align}
\delta_1\nu=\frac{1-\lambda^{-1}}{(1+\kappa^{-2})}t_1t_2, \qquad \delta_2\nu^{-1}=\frac{1-\lambda}{(1+\kappa^{2})}t_1t_2.   \label{eq:mu-and-inverse}
\end{align}
Their product yields $\delta_1\delta_2t_3^2=(2-\lambda-\lambda^{-1})t_1^2t_2^2$, so that
\begin{align}
t_3^2=\frac{(2-\lambda-\lambda^{-1})t_1^2t_2^2}{(t_1^2-2-\lambda-\lambda^{-1})(t_2^2-2-\lambda-\lambda^{-1})}.  \label{eq:case-3}
\end{align}

Set
$$\theta=\frac{(\lambda-1)(\kappa-\kappa^{-1})}{(\lambda+1)t_3}\ne 0,\pm\frac{\kappa-\kappa^{-1}}{t_3}.$$
Then
$$\lambda=\frac{\kappa-\kappa^{-1}+\theta t_3}{\kappa-\kappa^{-1}-\theta t_3}.$$
As is not difficult to verify, (\ref{eq:case-3}) becomes
$$t_3^2=\frac{4\theta^2(4+(\theta^2-1)t_3^2)t_1^2t_2^2t_3^2}{(\theta^2t_1^2t_3^2-(t_1^2-4)(t_3^2-4))(\theta^2t_2^2t_3^2-(t_2^2-4)(t_3^2-4))},$$
which, due to $t_3\ne 0,\pm 2$, is equivalent to
\begin{align*}
\theta^4+(4t_1^{-2}+4t_2^{-2}+4t_3^{-2}-2)\theta^2+(1-4t_1^{-2})(1-4t_2^{-2})(1-4t_3^{-2})=0.
\end{align*}

Amazingly, the $t_{ij}$'s and $t_{123}$ are determined in an elegant way:
\begin{align}
t_{ij}=\frac{\theta+1}{2}t_it_j, \quad  1\le i<j\le 3;  \qquad  t_{123}=\Big(\frac{\theta+1}{2}\Big)^2t_1t_2t_3.  \label{eq:trace-expression}
\end{align}
Indeed, we can easily compute
\begin{align*}
t_{13}&={\rm tr}\big(\mathbf{h}_{t_1}^{\lambda}(1)\mathbf{d}(\kappa)\big)=\frac{\lambda\kappa+\kappa^{-1}}{\lambda+1}t_1
=\frac{\theta+1}{2}t_1t_3, \\
t_{23}&={\rm tr}\big(\mathbf{h}_{t_2}^{\lambda}(\nu)\mathbf{d}(\kappa)\big)=\frac{\lambda\kappa+\kappa^{-1}}{\lambda+1}t_2
=\frac{\theta+1}{2}t_2t_3,   \\
t_{12}&={\rm tr}\big(\mathbf{h}_{t_1}^{\lambda}(1)\mathbf{h}_{t_2}^{\lambda}(\nu)\big)
=\frac{(\lambda^2+1)t_1t_2+\lambda\delta_2\nu^{-1}+\lambda\delta_1\nu}{(\lambda+1)^2},  \\
t_{123}&={\rm tr}\big(\mathbf{h}_{t_1}^{\lambda}(1)\mathbf{h}_{t_2}^{\lambda}(\nu)\mathbf{d}(\kappa)\big)
=\frac{\kappa(\lambda^2t_1t_2+\lambda\delta_2\nu^{-1})+\kappa^{-1}(t_1t_2+\lambda\delta_1\nu)}{(\lambda+1)^2}.
\end{align*}
Using (\ref{eq:mu-and-inverse}), we further obtain
\begin{align*}
\frac{t_{12}}{t_1t_2}&=\frac{1}{(\lambda+1)^2}\Big(\lambda^2+1+\frac{\lambda-\lambda^2}{1+\kappa^2}+\frac{\lambda-1}{1+\kappa^{-2}}\Big)
=\frac{\lambda\kappa+\kappa^{-1}}{(\lambda+1)t_3}=\frac{\theta+1}{2},   \\
\frac{t_{123}}{t_1t_2t_3}&=\frac{1}{(\lambda+1)^2t_3}\Big(\kappa\lambda^2+\frac{\kappa(\lambda-\lambda^2)}{1+\kappa^2}
+\kappa^{-1}+\frac{\kappa^{-1}(\lambda-1)}{1+\kappa^{-2}}\Big)  \\
&=\frac{(\lambda\kappa+\kappa^{-1})^2}{(\lambda+1)^2t_3^2}=\Big(\frac{\theta+1}{2}\Big)^2.
\end{align*}

\subsubsection{$\lambda_3=-1$} \label{sec:lambda3-special}

Now that $\mathbf{g}_3\ne-\mathbf{e}$, up to conjugacy we may assume $\mathbf{g}_3=-\mathbf{p}$, so $\mathbf{x}_3=\epsilon\mathbf{p}(\theta/2)$ for some $\epsilon\in\{\pm 1\}, \theta\ne 0$.

\begin{lemma}
$t_1t_2\ne 0$.
\end{lemma}
\begin{proof}
If $t_1=0$, then since $\mathbf{x}_1\cdot\overline{\mathbf{x}_2}\lrcorner\overline{\mathbf{x}_1}=\mathbf{g}_3=-\mathbf{p}$, by Lemma \ref{lem:key} (c') we have $\mathbf{x}_1\in \mathcal{T}_+$, so that $\mathbf{g}_2=[\mathbf{x}_3,\overline{\mathbf{x}_1}]\in\mathcal{T}_+$, which implies $\mathbf{x}_2\in{\rm Cen}(\mathbf{g}_2)\subset\mathcal{T}_+$. This contradicts the irreducibility.

Similarly, $t_2=0$ is neither possible.
\end{proof}

An application of Lemma \ref{lem:commutator} (c) to $\mathbf{a}_1=\mathbf{x}_1$ and $\mathbf{a}_2=\overline{\mathbf{x}_2}$ shows $\mathbf{x}_1=\mathbf{k}_{t_1}(\alpha)$, $\mathbf{x}_2=\mathbf{k}_{t_2}(\beta)$ for some $\alpha,\beta$; by conjugating via $\mathbf{p}(-\alpha(2t_1)^{-1})$, we may just assume $\mathbf{x}_1=\mathbf{k}_{t_1}(0)$. Hence the condition in Lemma \ref{lem:commutator} (c) reads
\begin{align}
t_1^{-1}t_2+t_1t_2^{-1}(\beta^2+1)=-t_1\beta.  \label{eq:case2-1}
\end{align}
Applying Lemma \ref{lem:commutator} (c) to $\mathbf{a}_1=\mathbf{x}_3\lrcorner\mathbf{x}_1=\mathbf{k}_{t_1}(t_1\theta)$ and $\mathbf{a}_2=\overline{\mathbf{x}_2}$, we see
\begin{align}
t_1^{-1}t_2(t_1^2\theta^2+1)+t_1t_2^{-1}(\beta^2+1)=2t_1\theta\beta+t_1t_2\theta-t_1\beta.  \label{eq:case2-2}
\end{align}

Summing (\ref{eq:case2-1}) and (\ref{eq:case2-2}) yields $t_1t_2\theta^2=2t_1\theta\beta+t_1t_2\theta$. Hence
\begin{align}
\beta=\frac{\theta-1}{2}t_2.   \label{eq:beta}
\end{align}
Then it is not difficult to verify that (\ref{eq:case2-1}) is equivalent to
$$\theta^2=1-4t_1^{-2}-4t_2^{-2}.$$

Note that
\begin{align*}
t_{13}&={\rm tr}(\epsilon\mathbf{p}(\theta/2)\mathbf{k}_{t_1}(0))=\epsilon(\theta+1)t_1, \\
t_{23}&={\rm tr}(\epsilon\mathbf{p}(\theta/2)\mathbf{k}_{t_2}(\beta))=\epsilon(\theta+1)t_2.
\end{align*}
By direct computation,
$$\mathbf{x}_1\mathbf{x}_2=\mathbf{k}_{t_1}(0)\mathbf{k}_{t_2}(\beta)
=\frac{1}{2}\left(\begin{array}{cc} t_1t_2-2t_1^{-1}t_2+t_1\beta &  \ast \\
4(t_1\beta+t_1t_2)  & t_1t_2-2t_1t_2^{-1}(1+\beta^2)-t_1\beta \end{array}\right),$$
where $\ast$ stands for some irrelevant number. So
\begin{align*}
t_{12}&={\rm tr}(\mathbf{x}_1\mathbf{x}_2)=t_1t_2-\frac{t_2}{t_1}-\frac{t_1}{t_2}(\beta^2+1)\stackrel{(\ref{eq:case2-1})}
=t_1t_2+t_1\beta\stackrel{(\ref{eq:beta})}=\frac{\theta+1}{2}t_1t_2, \\
t_{123}&={\rm tr}(\mathbf{x}_1\mathbf{x}_2\mathbf{x}_3)
=\epsilon\Big(t_1t_2-\frac{t_2}{t_1}-\frac{t_1}{t_2}(\beta^2+1)+\theta(t_1\beta+t_1t_2)\Big)=\frac{\epsilon(\theta+1)^2}{2}t_1t_2.
\end{align*}
Hence (remembering $t_3=2\epsilon$) we obtain the same expressions as (\ref{eq:trace-expression}):
$$t_{ij}=\frac{\theta+1}{2}t_it_j, \quad 1\le i<j\le 3; \qquad t_{123}=\Big(\frac{\theta+1}{2}\Big)^2t_1t_2t_3.$$

\begin{remark}
\rm If $t_1=t_2=t_3=2$, then $\epsilon=1$, and $\theta=\varepsilon\sqrt{-1}$, $\beta=\varepsilon\sqrt{-1}-1$ for some $\varepsilon\in\{\pm1\}$,
so
$$\mathbf{x}_1=\left(\begin{array}{cc} 1 & 0 \\ 4 & 1 \end{array}\right),  \quad
\mathbf{x}_2=\left(\begin{array}{cc} \varepsilon\sqrt{-1} & \varepsilon\sqrt{-1}/2 \\ 4 & 2-\varepsilon\sqrt{-1} \end{array}\right), \quad
\mathbf{x}_3=\left(\begin{array}{cc} 1 & \varepsilon\sqrt{-1}/2 \\ 0 & 1 \end{array}\right).$$
Conjugated by $\mathbf{d}(\sqrt{2})$, they become
$$\mathbf{x}_1=\left(\begin{array}{cc} 1 & 0 \\ 2 & 1 \end{array}\right),  \quad
\mathbf{x}_2=\left(\begin{array}{cc} \varepsilon\sqrt{-1} & \varepsilon\sqrt{-1} \\ 2 & 2-\varepsilon\sqrt{-1} \end{array}\right), \quad
\mathbf{x}_3=\left(\begin{array}{cc} 1 & \varepsilon\sqrt{-1} \\ 0 & 1 \end{array}\right).$$
This recovers the main result of \cite{CE21}, where it was shown that up to conjugacy, there are exactly two discrete faithful representations into ${\rm PSL}(2,\mathbb{C})$ sending meridians to parabolic elements.
One may also refer to \cite{Ma06}.
\end{remark}

\subsubsection{Unifying the cases $\lambda_3\ne -1$ and $\lambda_3=-1$}

The result of Section \ref{sec:lambda3-special} can be incorporated into Section \ref{sec:lambda3-general} by allowing $t_3=\pm 2$.
Indeed, for each $(t_1,t_2,t_3,\theta)$ with 
\begin{align}
\theta^4+(4t_1^{-2}+4t_2^{-2}+4t_3^{-2}-2)\theta^2+(1-4t_1^{-2})(1-4t_2^{-2})(1-4t_3^{-2})=0,   \label{eq:char-var-main}
\end{align}
the character determined by
\begin{align}
t_{ij}=\frac{\theta+1}{2}t_it_j, \quad 1\le i<j\le 3; \qquad  t_{123}=\Big(\frac{\theta+1}{2}\Big)^2t_1t_2t_3   \label{eq:expression'}
\end{align}
can be realized by the representation given by
\begin{align}
\mathbf{x}_1&=\frac{1}{2\theta t_1t_3}\left(\begin{array}{cc} \theta t_1^2t_3+(\kappa-\kappa^{-1})(t_1^2-4) &  t_1^2-4  \\
 \theta^2t_1^2t_3^2-(t_1^2-4)(t_3^2-4) & \theta t_1^2t_3-(\kappa-\kappa^{-1})(t_1^2-4) \end{array}\right),  \label{eq:realization-1} \\
\mathbf{x}_2&=\frac{t_2}{t_3}\left(\begin{array}{cc}  \kappa+(\theta-1)t_3/2 &  \star \\  2+\kappa^{-1}(\theta-1)t_3 & \kappa^{-1}+(1-\theta)t_3/2 \end{array}\right), \label{eq:realization-2}  \\
\mathbf{x}_3&=\left(\begin{array}{cc} \kappa & 1 \\ 0 & \kappa^{-1} \end{array} \right),  \label{eq:realization-3}
\end{align}
where $\kappa+\kappa^{-1}=t_3$, and $\star$ is the number determined by the condition $\det(\mathbf{x}_2)=1$ (its expression is too complicated, so we do not write it down). Thus, the cases $\lambda_3\ne -1$ and $\lambda_3=-1$ can be unified.

The geometry turns out to be extremely simple: rewriting (\ref{eq:char-var-main}) as
$$4\big((4t_1^{-2}-1)(4t_2^{-2}-1)-\theta^2\big)=(\theta^2+4t_1^{-2}-1)(\theta^2+4t_2^{-2}-1)t_3^2,$$
we see
$$\theta^2\ne 1-4t_1^{-2},\ 1-4t_2^{-2}, \ (1-4t_1^{-2})(1-4t_2^{-2}),$$
so that
$$\frac{4}{t_3^2}=\frac{(\theta^2+4t_1^{-2}-1)(\theta^2+4t_2^{-2}-1)}{(4t_1^{-2}-1)(4t_2^{-2}-1)-\theta^2}.$$
Thus, the locus of (\ref{eq:char-var-main}) is a $2$-fold regular cover over
$$\big\{(t_1,t_2,\theta)\in(\mathbb{C}^\ast)^3\colon \theta^2\ne 1-4t_1^{-2},\ 1-4t_2^{-2}, \ (1-4t_1^{-2})(1-4t_2^{-2})\big\}.$$
Due to the symmetric nature, the projection to the $(t_1,t_3,\theta)$-space or $(t_2,t_3,\theta)$-space also defines a $2$-fold regular cover.

\subsection{The result}

It is known (see \cite{Go09} Section 5) that $\mathcal{X}(F_3)$ is isomorphic to
\begin{align}
\{\vec{t}=(t_1,t_2,t_3,t_{12},t_{13},t_{23}, t_{123})\colon t_{123}^2-\nu_1t_{123}+\nu_0=0\}\subset\mathbb{C}^7,   \label{eq:F3}
\end{align}
where
\begin{align*}
\nu_0&=t_1^2+t_2^2+t_3^2+t_{12}^2+t_{13}^2+t_{23}^2-t_1t_2t_{12}-t_1t_3t_{13}-t_2t_3t_{23}+t_{12}t_{13}t_{23}-4, \\
\nu_1&=t_1t_{23}+t_2t_{13}+t_3t_{12}-t_1t_2t_3.
\end{align*}
It means that the ${\rm SL}(2,\mathbb{C})$-character variety of any 3-generator group can be embedded into the hypersurface given by (\ref{eq:F3}).

For $i\ne j$, let
$$f_{i,j}=t_i^2+t_j^2+t_{ij}^2-t_it_jt_{ij}-2.$$
Then $f_{i,j}={\rm tr}([\mathbf{x}_i,\mathbf{x}_j])={\rm tr}([\mathbf{x}_i,\overline{\mathbf{x}_j}])$, as seen from
\begin{align*}
{\rm tr}([\mathbf{x}_i,\mathbf{x}_j])&={\rm tr}(\mathbf{x}_i\mathbf{x}_j(\mathbf{x}_j\mathbf{x}_i)^{-1})
={\rm tr}\big(\mathbf{x}_i\mathbf{x}_j(t_{ij}\mathbf{e}-\mathbf{x}_j\mathbf{x}_i)\big)  \\
&=t_{ij}^2-{\rm tr}(\mathbf{x}_i^2\mathbf{x}_j^2)=t_{ij}^2-{\rm tr}\big((t_i\mathbf{x}_i-\mathbf{e})(t_j\mathbf{x}_j-\mathbf{e})\big)=f_{i,j},  \\
{\rm tr}([\mathbf{x}_i,\overline{\mathbf{x}_j}])&={\rm tr}(\mathbf{x}_i\overline{\mathbf{x}_j}\overline{\mathbf{x}_i}\mathbf{x}_j)
={\rm tr}(\mathbf{x}_j\mathbf{x}_i\overline{\mathbf{x}_j}\overline{\mathbf{x}_i})={\rm tr}([\mathbf{x}_j,\mathbf{x}_i]),
\end{align*}
and that $f_{i,j}$ is symmetric in $i,j$.
By Lemma \ref{lem:commutator} (e), $f_{i,j}=2$ if and only if $\mathbf{x}_i$ and $\mathbf{x}_j$ have a common eigenvector.

\begin{thm}
The irreducible ${\rm SL}(2,\mathbb{C})$-character variety of the Borromean link can be decomposed as
$$\mathcal{X}^{\rm irr}(B)=\big(\cup_{i=1}^3\mathcal{X}^{+}_{1,i}\big)\cup\big(\cup_{i=1}^3\mathcal{X}^{-}_{1,i}\big)
\cup\big(\cup_{i=1}^3\mathcal{X}_{2,i}\big)\cup\big(\cup_{i=1}^3\mathcal{X}_{3,i}\big)\cup\mathcal{X}_{4},$$
where
\begin{align*}
\mathcal{X}^+_{1,i}&\cong\big\{\vec{t}\colon t_i=2,\ t_{i,i\pm1}=t_{i\pm1},\ t_{123}=t_{i-1,i+1}, \ f_{i-1,i+1}\ne 2\big\},  \\
\mathcal{X}^-_{1,i}&\cong\big\{\vec{t}\colon t_i=-2,\ t_{i,i\pm1}=-t_{i\pm1},\ t_{123}=-t_{i-1,i+1}, \ f_{i-1,i+1}\ne 2\big\},  \\
\mathcal{X}_{2,i}&\cong\big\{\vec{t}\colon t_i=t_{i,i\pm1}=t_{123}=0,\ f_{i-1,i+1}=2, \ t_{i\pm1}^2\ne 4\big\},   \\
\mathcal{X}_{3,i}&\cong\big\{\vec{t}\colon t_{i\pm1}=t_{i-1,i+1}=0,\ t_it_{123}=t_{i,i-1}t_{i,i+1}, \
t_{i,i-1}^2+t_{i,i+1}^2+t_{123}^2=4-t_i^2\ne 0\big\},  \\
\mathcal{X}_4&\cong\Big\{(t_1,t_2,t_3,\theta)\in(\mathbb{C}^\ast)^4\colon \theta^4-\Big(2-4{\sum}_{i=1}^3t_i^{-2}\Big)\theta^2+{\prod}_{i=1}^3(1-4t_i^{-2})=0\Big\}.
\end{align*}
\end{thm}

Some supplements are in order.
\begin{enumerate}
  \item Characters in $\mathcal{X}^{+}_{1,i}$ (resp. $\mathcal{X}^{-}_{1,i}$) are realized by representations such that
        $\mathbf{x}_i=\mathbf{e}$ (resp. $\mathbf{x}_i=-\mathbf{e}$) and there is no constraint on $\mathbf{x}_{i-1},\mathbf{x}_{i+1}$, except that $\mathbf{x}_{i-1}$, $\mathbf{x}_{i+1}$ have no common eigenvector. 
  \item Characters in $\mathcal{X}_{2,i}$, $\mathcal{X}_{3,i}$ are respectively realized by the representations given in Lemma \ref{lem:g=e}
        and Lemma \ref{lem:g=-e}.
  \item Characters in $\mathcal{X}_4$ are realized by representations given by (\ref{eq:realization-1})--(\ref{eq:realization-3}).
        With (\ref{eq:char-var-main}), (\ref{eq:expression'}), it is ensured that for all $i$,
        $$\lambda_i+\lambda_i^{-1}=f_{i-1,i+1}=\frac{\theta^2-1}{4}t_{i-1}^2t_{i+1}^2+t_{i-1}^2+t_{i+1}^2-2\ne 2,$$
        so that $\mathbf{g}_i\ne\mathbf{e}$.
        Moreover, $t_1t_2t_3\ne 0$ ensures $\mathbf{g}_i\ne-\mathbf{e}$ for all $i$.

        Since $\mathcal{X}_4$ is irreducible and contains the character of a lift of the {\it holonomy representation}
        $\pi(B)\to{\rm PSL}(2,\mathbb{C})$ (i.e. the one defining the hyperbolic structure of $S^3\setminus B$), we see that $\mathcal{X}_4$ is the canonical component.
        It can be exhibited as a regular $2$-cover over the complement of three hypersurfaces in $(\mathbb{C}^\ast)^3$, in three ways.
  \item The nonempty intersections among the various parts are: for each $i$,
        \begin{align*}
        \mathcal{X}_{3,i}\cap\mathcal{X}_{2,i+1}&=\big\{\vec{t}\colon t_{i\pm1}=t_{i-1,i+1}=t_{i,i+1}=t_{123}=0,\ t_{i,i-1}^2=4-t_i^2\ne 0\big\},\\
        \mathcal{X}_{3,i}\cap\mathcal{X}_{2,i-1}&=\big\{\vec{t}\colon t_{i\pm1}=t_{i-1,i+1}=t_{i,i-1}=t_{123}=0,\ t_{i,i+1}^2=4-t_i^2\ne 0\big\},\\
        \mathcal{X}_{3,i-1}\cap\mathcal{X}_{3,i+1}&=\big\{\vec{t}\colon t_1=t_2=t_3=t_{i,i\pm1}=0,\ t_{123}^2+t_{i-1,i+1}^2=4\big\}.
        \end{align*}
\end{enumerate}

\section{Twisted Alexander polynomials}

We briefly recall the definition of twisted Alexander polynomial, with small modifications. One may refer to \cite{Wa94} for more detail.

For a ring $R$, let $\mathcal{M}_n(R)$ denote the ring of $n\times n$ matrices over $R$.

Let $L=K_1\sqcup\ldots\sqcup K_m$ be an oriented $m$-component link, and let
$$\pi:=\pi(L)\cong\langle x_1,\ldots,x_\ell\mid r_1,\ldots,r_{\ell-1}\rangle$$
be a presentation such that each $x_j$ comes from an arc of $K_{\sigma(j)},\sigma(j)\in\{1,\ldots,m\}$.
Let $F_\ell=\langle x_1,\ldots,x_\ell\mid -\rangle$, the free group generated by $x_1,\ldots,x_\ell$.
Let $M$ be the $(\ell-1)\times\ell$ matrix whose $(i,j)$-entry is the image of $\partial r_i/\partial x_j$ (the Fox derivative) under the ring homomorphism $$q:\mathbb{Z}[F_\ell]\to\mathbb{Z}[\pi]$$
induced by the canonical map $F_\ell\to\pi$, and let $M_v\in\mathcal{M}_{\ell-1}(\mathbb{Z}[\pi])$ be the matrix obtained from deleting the $v$-th column of $M$.
Let
$$\mathfrak{a}:\pi\to \mathbb{Z}^{\oplus m}=\langle s_1\rangle\oplus\cdots\oplus\langle s_m\rangle$$
denote the abelianization map, which sends $x_j$ to $s_{\sigma(j)}$.

For a representation $\rho:\pi\to{\rm SL}(2,\mathbb{C})$, the composite
$$\pi\stackrel{\mathfrak{a}\times\rho}\longrightarrow\mathbb{Z}^{\oplus m}\times {\rm SL}(2,\mathbb{C})\hookrightarrow \mathbb{Z}[s_1^{\pm1},\ldots,s_m^{\pm1}]\times\mathcal{M}_2(\mathbb{C})\to\mathcal{M}_2(\mathbb{C}[s_1^{\pm1},\ldots,s_m^{\pm1}])$$
can be extended by linearity to a ring homomorphism
\begin{align*}
\Phi:\mathbb{Z}[\pi]\to\mathcal{M}_2(\mathbb{C}[s_1^{\pm1},\ldots,s_m^{\pm1}]).
\end{align*}
The {\it twisted Alexander polynomial} of $L$ associated to $\rho$ is defined to be
$$\Delta_{L}^{\rho}=\Delta_{L}^{\rho}(s_1,\ldots,s_m)\doteq\frac{\det\Phi(M_v)}{\det\Phi(1-x_v)}\in \mathbb{C}(s_1^{\pm1},\ldots,s_m^{\pm 1}),$$
where $\Phi(M_v)\in\mathcal{M}_{2(\ell-1)}(\mathbb{C}[s_1^{\pm1},\ldots,s_m^{\pm1}])$ is the big matrix obtained from $M_v$ by replacing each entry with its image under $\Phi$, and $\doteq$ means an equality up to multiplication by $s_1^{k_1}\cdots s_m^{k_m}$ for $k_1,\ldots,k_m\in\mathbb{Z}$.
This does not depend on $v$, neither on the presentation as long as it is {\it strongly Tietze equivalent} to a Wirtinger presentation.
It is known that $\Delta_{L}^{\rho}$ is a Laurent polynomial when $m\ge 2$; see \cite{Wa94} Proposition 9.

Actually, $\Delta_L^{\rho}$ only depends on the conjugacy class of $\rho$, so we also denote it by $\Delta_L^{[\rho]}$. Here $[\rho]$ stands for the conjugacy class of $\rho$, which can be identified with $\chi_\rho$ if $\rho$ is irreducible.
This manifests TAP as a function on $\mathcal{X}^{\rm irr}(L)$.

\bigskip

Now let $L=B$. The presentation (\ref{eq:presentation}) is strongly Tietze equivalent to (\ref{eq:Wirtinger}), so (\ref{eq:presentation}) can be used to define $\Delta_B^{[\rho]}$. 

\begin{thm}\label{thm:TAP}
The twisted Alexander polynomial of the Borromean link associated to an irreducible representation is given by
$$\Delta^{[\rho]}_B\doteq \prod_{j=1}^3(s_j+s_j^{-1}-t_j)
+\begin{cases}
0, & [\rho]\in\mathcal{X}^{\pm}_{1,i}\ \text{or\ }[\rho]\in\mathcal{X}_{2,i},  \\
4t_{123}+t_1t_2t_3, &[\rho]\in\mathcal{X}_{3,i}, \\
\theta^2t_1t_2t_3, &[\rho]\in\mathcal{X}_4. \end{cases}$$
\end{thm}

\begin{remark}
\rm The {\it hyperbolic torsion conjecture} proposed in \cite{DFJ12} states that for a hyperbolic knot $K$, the TAP associated to a lift $\rho_0$ of the holonomy representation detects the genus and fibredness of $K$. It was generalized to hyperbolic links in \cite{MT17}.
According to \cite{MT17} Remark 3.4, the conjecture for a $m$-component alternating link $L$ is equivalent to $\deg\Delta_L^{\rho_0}=4g(L)+2(m-2).$

It is known that the Borromean link is fibred with genus 1. Thus, we have confirmed the hyperbolic torsion conjecture for the Borromean link.
\end{remark}

The remaining part is devoted to proving Theorem \ref{thm:TAP}.

To keep the expressions compact, for $\alpha\in\mathbb{Z}[F_3]$ we denote $q(\alpha)\in\mathbb{Z}[\pi]$ also by $\alpha$, where $\pi=\pi(B)$.
With the presentation (\ref{eq:presentation}),
$r_1=[x_2,[x_3,\overline{x_1}]]$, $r_2=[x_3,[x_1,\overline{x_2}]]$.
Remembering the property of Fox derivative (cf. \cite{GTM175} Page 117) and using that when $r_i=f\overline{g}$ with $f,g\in F_3$,
$$\frac{\partial r_i}{\partial x_j}=\frac{\partial}{\partial x_j}(f\overline{g})=\frac{\partial}{\partial x_j}f-f\overline{g}\frac{\partial}{\partial x_j}g=\frac{\partial}{\partial x_j}(f-g) \quad \text{in\ } \mathbb{Z}[\pi],$$
we can compute
\begin{align*}
\frac{\partial r_1}{\partial x_1}&=\frac{\partial}{\partial x_1}\big((x_2x_3\overline{x_1}\overline{x_3}x_1)(x_3\overline{x_1}\overline{x_3}x_1x_2)^{-1}\big)
=\frac{\partial}{\partial x_1}(x_2x_3\overline{x_1}\overline{x_3}x_1-x_3\overline{x_1}\overline{x_3}x_1x_2) \\
&=x_2x_3(-\overline{x_1})+x_2x_3\overline{x_1}x_3-x_3(-\overline{x_1})-x_3\overline{x_1}\overline{x_3}=(1-x_2)x_3\overline{x_1}(1-\overline{x_3}),  \\
\frac{\partial r_1}{\partial x_2}&=\frac{\partial}{\partial x_2}(x_2x_3\overline{x_1}\overline{x_3}x_1-x_3\overline{x_1}\overline{x_3}x_1x_2)=1-[x_3,\overline{x_1}],  \\
\frac{\partial r_2}{\partial x_1}&=\frac{\partial}{\partial x_1}(x_3x_1\overline{x_2}\overline{x_1}x_2-x_1\overline{x_2}\overline{x_1}x_2x_3) \\
&=x_3+x_3x_1\overline{x_2}(-\overline{x_1})-1-x_1\overline{x_2}(-\overline{x_1})=(1-x_3)(x_1\lrcorner\overline{x_2}-1),  \\
\frac{\partial r_2}{\partial x_2}&=x_3x_1(-\overline{x_2})+x_3x_1\overline{x_2}\overline{x_1}-x_1(-\overline{x_2})-x_1\overline{x_2}\overline{x_1}
=(1-x_3)x_1\overline{x_2}(1-\overline{x_1}).
\end{align*}
Hence
\begin{align}
M_3=\left(\begin{array}{cc} (1-x_2)x_3\overline{x_1}(1-\overline{x_3}) & 1-[x_3,\overline{x_1}] \\
(1-x_3)(x_1\lrcorner\overline{x_2}-1) & (1-x_3)x_1\overline{x_2}(1-\overline{x_1}) \end{array}\right).   \label{eq:M3}
\end{align}

For a $2\times 2$ matrix $\mathbf{u}$ over some commutative unitary ring, denote its adjoint by $\mathbf{u}^\ast$, so
$$\left(\begin{array}{cc} a & b \\ c & d \end{array}\right)^\ast=\left(\begin{array}{cc} d & -b \\ -c & a \end{array}\right).$$
We shall use the following facts:
$$\mathbf{u}\mathbf{u}^\ast=\det(\mathbf{u})\mathbf{e},  \quad \mathbf{u}+\mathbf{u}^\ast={\rm tr}(\mathbf{u})\mathbf{e},  \quad
(\mathbf{u}+\mathbf{v})^\ast=\mathbf{u}^\ast+\mathbf{v}^\ast, \quad (\mathbf{u}\mathbf{v})^\ast=\mathbf{v}^\ast\mathbf{u}^\ast.$$

\begin{lemma}\label{lem:matrix}
If $\mathbf{x}\in{\rm SL}(2,\mathbb{C})$ with ${\rm tr}(\mathbf{x})=t$, then
$$\det(s\mathbf{x}-\mathbf{e})=\det(s\mathbf{e}-\mathbf{x})=s^2-ts+1.$$
\end{lemma}

\begin{proof}
We have
\begin{align*}
\det(s\mathbf{x}-\mathbf{e})\cdot\mathbf{e}&=(s\mathbf{x}-\mathbf{e})(s\mathbf{x}-\mathbf{e})^\ast
=(s\mathbf{x}-\mathbf{e})(s\mathbf{x}^\ast-\mathbf{e})   \\
&=s^2\mathbf{x}\mathbf{x}^\ast-s(\mathbf{x}+\mathbf{x}^\ast)+\mathbf{e}=(s^2-ts+1)\mathbf{e}.
\end{align*}
Similarly, $\det(s\mathbf{e}-\mathbf{x})=s^2-ts+1$.
\end{proof}

Suppose $\rho(x_i)=\mathbf{x}_i$, $i=1,2,3$. By (\ref{eq:M3}),
\begin{align}
\Delta^{[\rho]}_B&\doteq\frac{\det\Phi(M_3)}{\det\Phi(1-x_3)}
\doteq\det\left(\begin{array}{cc} (\mathbf{e}-s_2\mathbf{x}_2)\mathbf{x}_3\overline{\mathbf{x}_1}(s_3\mathbf{e}-\overline{\mathbf{x}_3}) & s_1(\mathbf{e}-[\mathbf{x}_3,\overline{\mathbf{x}_1}])  \\
\mathbf{x}_1\lrcorner\overline{\mathbf{x}_2}-s_2\mathbf{e} & \mathbf{x}_1\overline{\mathbf{x}_2}(s_1\mathbf{e}-\overline{\mathbf{x}_1}) \end{array}\right)    \nonumber \\
&\doteq\det\left(\begin{array}{cc} (\mathbf{e}-s_2\mathbf{x}_2)\mathbf{x}_3\overline{\mathbf{x}_1}(s_3\mathbf{e}-\overline{\mathbf{x}_3}) & s_1(\mathbf{e}-[\mathbf{x}_3,\overline{\mathbf{x}_1}])  \\ \overline{\mathbf{x}_1}-s_2\mathbf{x}_2\overline{\mathbf{x}_1} & s_1\mathbf{e}-\overline{\mathbf{x}_1} \end{array}\right)  \nonumber  \\
&\doteq\det\left(\begin{array}{cc} s_3(\mathbf{e}-s_2\mathbf{x}_2)\mathbf{x}_3+[\mathbf{x}_3,\overline{\mathbf{x}_1}](s_2\mathbf{x}_2-\mathbf{e})  & s_1(\mathbf{e}-[\mathbf{x}_3,\overline{\mathbf{x}_1}])  \\ \mathbf{e}-s_2\mathbf{x}_2 & s_1\mathbf{e}-\overline{\mathbf{x}_1} \end{array}\right)
\nonumber \\
&=\det\left(\begin{array}{cc} s_3(\mathbf{e}-s_2\mathbf{x}_2)\mathbf{x}_3  & s_1\mathbf{e}-\mathbf{x}_3\lrcorner\overline{\mathbf{x}_1}\\
\mathbf{e}-s_2\mathbf{x}_2 & s_1\mathbf{e}-\overline{\mathbf{x}_1} \end{array}\right) \nonumber \\
&\doteq\det\big(s_3(s_1\mathbf{e}-\mathbf{x}_3\lrcorner\mathbf{x}_1)(\mathbf{e}-s_2\mathbf{x}_2)\mathbf{x}_3
+(s_1\mathbf{e}-\mathbf{x}_1)(s_2\mathbf{x}_2-\mathbf{e})\big). \label{eq:TAP}
\end{align}
Here the third line is obtained by multiplying the first column (in the matrix) by $\mathbf{x}_1$ on the right,
the fourth line results from adding $[\mathbf{x}_3,\overline{\mathbf{x}_1}]$ times the second row to the first row,
and the last line is deduced as follows.
Let
$$\mathbf{u}_1=s_3(\mathbf{e}-s_2\mathbf{x}_2)\mathbf{x}_3, \ \ \ \mathbf{u}_2=s_1\mathbf{e}-\mathbf{x}_3\lrcorner\overline{\mathbf{x}_1}, \ \ \
\mathbf{u}_3=\mathbf{e}-s_2\mathbf{x}_2, \ \ \ \mathbf{u}_4=s_1\mathbf{e}-\overline{\mathbf{x}_1},$$
regarded as elements of ${\rm GL}\big(2,\mathbb{C}(s_1^{\pm1},s_2^{\pm1},s_3^{\pm 1})\big)$,
then
\begin{align*}
\Delta^{[\rho]}_B&\doteq
\det\left(\begin{array}{cc}  \mathbf{u}_1 & \mathbf{u}_2 \\ \mathbf{u}_3 & \mathbf{u}_4 \end{array}\right)
=\det(\mathbf{u}_2)\det(\mathbf{u}_4)\det\left(\begin{array}{cc}  \overline{\mathbf{u}_2}\mathbf{u}_1 & \mathbf{e} \\
\overline{\mathbf{u}_4}\mathbf{u}_3 & \mathbf{e} \end{array}\right)  \\
&=\det(\mathbf{u}_2)\det(\mathbf{u}_4)\det(\overline{\mathbf{u}_2}\mathbf{u}_1-\overline{\mathbf{u}_4}\mathbf{u}_3)
=\det(\mathbf{u}_2^\ast\mathbf{u}_1-\mathbf{u}_4^\ast\mathbf{u}_3),
\end{align*}
using $\det(\mathbf{u}_2)=\det(\mathbf{u}_4)$. Hence (\ref{eq:TAP}) holds.

\begin{lemma}
If $\mathbf{x}_k\mathbf{x}_j=\mathbf{x}_j\mathbf{x}_k$ for some $j\ne k$, then
\begin{align}
\Delta^{[\rho]}_L\doteq\det((s_1\mathbf{e}-\mathbf{x}_1)(\mathbf{e}-s_2\mathbf{x}_2)(s_3\mathbf{x}_3-\mathbf{e}))={\prod}_{j=1}^3(s_j^2-t_js_j+1).
\end{align}
This happens when $[\rho]\in\mathcal{X}^{\pm}_{1,i}$ or $[\rho]\in\mathcal{X}_{2,i}$.
\end{lemma}

\begin{proof}
Let
$$\mathbf{w}=s_3(s_1\mathbf{e}-\mathbf{x}_3\lrcorner\mathbf{x}_1)(\mathbf{e}-s_2\mathbf{x}_2)\mathbf{x}_3
+(s_1\mathbf{e}-\mathbf{x}_1)(s_2\mathbf{x}_2-\mathbf{e}).$$

When $\mathbf{x}_1\mathbf{x}_3=\mathbf{x}_3\mathbf{x}_1$, we have $\mathbf{x}_3\lrcorner\mathbf{x}_1=\mathbf{x}_1$, so the assertion is obvious.

When $\mathbf{x}_2\mathbf{x}_3=\mathbf{x}_3\mathbf{x}_2$, to verify the assertion, just notice
\begin{align*}
\mathbf{w}&=\big(s_3(s_1\mathbf{e}-\mathbf{x}_3\lrcorner\mathbf{x}_1)\mathbf{x}_3-(s_1\mathbf{e}-\mathbf{x}_1)\big)(\mathbf{e}-s_2\mathbf{x}_2)  \\
&=(s_3\mathbf{x}_3-\mathbf{e})(s_1\mathbf{e}-\mathbf{x}_1)(\mathbf{e}-s_2\mathbf{x}_2).
\end{align*}

When $\mathbf{x}_1\mathbf{x}_2=\mathbf{x}_2\mathbf{x}_1$, from $\mathbf{x}_2[\mathbf{x}_3,\overline{\mathbf{x}_1}]=[\mathbf{x}_3,\overline{\mathbf{x}_1}]\mathbf{x}_2$ we can deduce that $\mathbf{x}_2$ commutes with $\mathbf{x}_3\lrcorner\overline{\mathbf{x}_1}$. Then
\begin{align*}
\mathbf{w}&=(\mathbf{e}-s_2\mathbf{x}_2) \big(s_3(s_1\mathbf{e}-\mathbf{x}_3\lrcorner\mathbf{x}_1)\mathbf{x}_3-(s_1\mathbf{e}-\mathbf{x}_1)\big) \\
&=(\mathbf{e}-s_2\mathbf{x}_2)(s_3\mathbf{x}_3-\mathbf{e})(s_1\mathbf{e}-\mathbf{x}_1).
\end{align*}
Hence the assertion holds.
\end{proof}

Go on to deal with (\ref{eq:TAP}) for $[\rho]\in\mathcal{X}_{3,i}$ or $[\rho]\in\mathcal{X}_4$.
Let
\begin{align*}
\mathbf{u}&=(s_1\mathbf{e}-\mathbf{x}_3\lrcorner\mathbf{x}_1)(\mathbf{e}-s_2\mathbf{x}_2)\mathbf{x}_3
=-s_1s_2\mathbf{x}_2\mathbf{x}_3+s_1\mathbf{x}_3+s_2\mathbf{x}_3\mathbf{x}_1\overline{\mathbf{x}_3}\mathbf{x}_2\mathbf{x}_3-\mathbf{x}_3\mathbf{x}_1,  \\
\mathbf{v}&=(s_1\mathbf{e}-\mathbf{x}_1)(s_2\mathbf{x}_2-\mathbf{e}))
=s_1s_2\mathbf{x}_2-s_1\mathbf{e}-s_2\mathbf{x}_1\mathbf{x}_2+\mathbf{x}_1.
\end{align*}
Then
\begin{align*}
\det(s_3\mathbf{u}+\mathbf{v})\cdot\mathbf{e}&=(s_3\mathbf{u}+\mathbf{v})(s_3\mathbf{u}^\ast+\mathbf{v}^\ast)
=s_3^2\mathbf{u}\mathbf{u}^\ast+\mathbf{v}\mathbf{v}^\ast+s_3(\mathbf{u}\mathbf{v}^\ast+\mathbf{v}\mathbf{u}^\ast)  \\
&=(s_3^2\det(\mathbf{u})+\det(\mathbf{v}))\cdot\mathbf{e}+s_3(\mathbf{u}\mathbf{v}^\ast+\mathbf{v}\mathbf{u}^\ast);
\end{align*}
taking trace and dividing by $2$, we obtain
\begin{align}
\det(s_3\mathbf{u}+\mathbf{v})=s_3^2\det(\mathbf{u})+\det(\mathbf{v})+s_3{\rm tr}(\mathbf{u}\mathbf{v}^\ast). \label{eq:det}
\end{align}

We have
\begin{align*}
{\rm tr}(\mathbf{u}\mathbf{v}^\ast)
&=-s_1^2s_2^2{\rm tr}(\mathbf{x}_2\lrcorner\mathbf{x}_3)+s_1^2s_2{\rm tr}(\mathbf{x}_3\overline{\mathbf{x}_2}+\mathbf{x}_2\mathbf{x}_3)
+s_1s_2^2\eta_1-s_1^2{\rm tr}(\mathbf{x}_3)  \\
&\ \ \ -s_1s_2\eta_2-s_2^2\eta_3+s_1{\rm tr}(\mathbf{x}_3\overline{\mathbf{x}_1}+\mathbf{x}_3\mathbf{x}_1)+s_2\eta_4-{\rm tr}(\mathbf{x}_3)  \\
&=-s_1^2s_2^2t_3+s_1^2s_2t_2t_3+s_1s_2^2\eta_1-s_1^2t_3-s_1s_2\eta_2-s_2^2\eta_3+s_1t_1t_3+s_2\eta_4-t_3,
\end{align*}
where
\begin{align*}
\eta_1&={\rm tr}(\mathbf{x}_2\mathbf{x}_3\overline{\mathbf{x}_2}\overline{\mathbf{x}_1}+\mathbf{x}_3\mathbf{x}_1[\overline{\mathbf{x}_3},\mathbf{x}_2]),   \\
\eta_2&={\rm tr}(\mathbf{x}_3\mathbf{x}_1[\overline{\mathbf{x}_3},\mathbf{x}_2]\overline{\mathbf{x}_1}),  \\
\eta_3&={\rm tr}(\mathbf{x}_3\mathbf{x}_1\overline{\mathbf{x}_3}\mathbf{x}_2\mathbf{x}_3\overline{\mathbf{x}_1}
        +\mathbf{x}_3\mathbf{x}_1\overline{\mathbf{x}_2}\overline{\mathbf{x}_1}),  \\
\xi&={\rm tr}(\mathbf{x}_2\mathbf{x}_3\overline{\mathbf{x}_1}+\mathbf{x}_3\overline{\mathbf{x}_2}\overline{\mathbf{x}_1}
+\mathbf{x}_3\mathbf{x}_1\overline{\mathbf{x}_3}\mathbf{x}_2\mathbf{x}_3+\mathbf{x}_3\mathbf{x}_1\overline{\mathbf{x}_2});
\end{align*}
the following has been used: ${\rm tr}(\mathbf{x}_3\overline{\mathbf{x}_2})={\rm tr}(\mathbf{x}_3(t_2\mathbf{e}-\mathbf{x}_2))=t_2t_3-t_{23}$,
and similarly, ${\rm tr}(\mathbf{x}_3\overline{\mathbf{x}_1})=t_1t_3-t_{13}$.

Using $\mathbf{x}_1[\overline{\mathbf{x}_3},\mathbf{x}_2]=[\overline{\mathbf{x}_3},\mathbf{x}_2]\mathbf{x}_1$,
$[\overline{\mathbf{x}_1},\mathbf{x}_3]\mathbf{x}_2=\mathbf{x}_2[\overline{\mathbf{x}_1},\mathbf{x}_3]$,
we can compute
\begin{align*}
\eta_1&={\rm tr}(\mathbf{x}_2\mathbf{x}_3\overline{\mathbf{x}_2}\overline{\mathbf{x}_1}
+\mathbf{x}_3[\overline{\mathbf{x}_3},\mathbf{x}_2]\mathbf{x}_1)
={\rm tr}(\mathbf{x}_2\mathbf{x}_3\overline{\mathbf{x}_2}(\overline{\mathbf{x}_1}+\mathbf{x}_1))
={\rm tr}(\mathbf{x}_2\mathbf{x}_3\overline{\mathbf{x}_2}\cdot t_1\mathbf{e})=t_1t_3, \\
\eta_2&={\rm tr}(\mathbf{x}_3[\overline{\mathbf{x}_3},\mathbf{x}_2])={\rm tr}(\mathbf{x}_2\mathbf{x}_3\overline{\mathbf{x}_2})=t_3, \\
\eta_3&={\rm tr}([\overline{\mathbf{x}_1},\mathbf{x}_3]\mathbf{x}_2\mathbf{x}_3)
+{\rm tr}(\mathbf{x}_3(\mathbf{x}_1\overline{\mathbf{x}_2}\overline{\mathbf{x}_1}))
={\rm tr}(\mathbf{x}_2[\overline{\mathbf{x}_1},\mathbf{x}_3]\mathbf{x}_3)
+{\rm tr}((\mathbf{x}_1\overline{\mathbf{x}_2}\overline{\mathbf{x}_1})\mathbf{x}_3)  \\
&={\rm tr}(\mathbf{x}_1\mathbf{x}_2\overline{\mathbf{x}_1}\mathbf{x}_3)
+{\rm tr}(\mathbf{x}_1\overline{\mathbf{x}_2}\overline{\mathbf{x}_1}\mathbf{x}_3)
={\rm tr}(\mathbf{x}_1(\mathbf{x}_2+\overline{\mathbf{x}_2})\overline{\mathbf{x}_1})\mathbf{x}_3)
={\rm tr}(t_2\mathbf{e}\cdot\mathbf{x}_3)=t_2t_3.
\end{align*}
Consequently,
$${\rm tr}(\mathbf{u}\mathbf{v}^\ast)=-s_1^2s_2^2t_3+s_1^2s_2t_2t_3+s_1s_2^2t_1t_3-s_1^2t_3-s_2^2t_3-s_1s_2\xi+s_1t_1t_3+s_2t_2t_3-t_3.$$
Therefore, by (\ref{eq:TAP}), (\ref{eq:det}),
\begin{align}
\Delta^{[\rho]}_B&\doteq s_3^2\det(\mathbf{u})+\det(\mathbf{v})+s_3{\rm tr}(\mathbf{u}\mathbf{v}^\ast)  \nonumber  \\
&=(s_3^2+1)(s_1^2+1-s_1t_1)(s_2^2+1-s_2t_2)  \nonumber  \\
& \ \ \ +s_3\big(s_1(s_2^2+1)t_1t_3+s_2(s_1^2+1)t_2t_3-(s_1^2+1)(s_2^2+1)t_3-s_1s_2\xi\big)  \nonumber \\
&\doteq{\prod}_{j=1}^3(s_j+s_j^{-1}-t_j)+t_1t_2t_3-\xi.   \label{eq:TAP-2}
\end{align}

The determination of $\xi$ proceeds as
\begin{align*}
{\rm tr}(\mathbf{x}_2\mathbf{x}_3\overline{\mathbf{x}_1})&={\rm tr}(\mathbf{x}_2\mathbf{x}_3(t_1\mathbf{e}-\mathbf{x}_1))=t_1t_{23}-t_{123}, \\
{\rm tr}(\mathbf{x}_3\overline{\mathbf{x}_2}\overline{\mathbf{x}_1})
&={\rm tr}(\mathbf{x}_3(t_{12}\mathbf{e}-\mathbf{x}_1\mathbf{x}_2))=t_3t_{12}-t_{123}, \\
{\rm tr}(\mathbf{x}_3\mathbf{x}_1\overline{\mathbf{x}_2})&={\rm tr}(\mathbf{x}_3\mathbf{x}_1(t_2\mathbf{e}-\mathbf{x}_2))=t_2t_{13}-t_{123},  \\
{\rm tr}(\mathbf{x}_3\mathbf{x}_1\overline{\mathbf{x}_3}\mathbf{x}_2\mathbf{x}_3)
&={\rm tr}\big(\mathbf{x}_3\mathbf{x}_1({\rm tr}(\mathbf{x}_2\overline{\mathbf{x}_3})\mathbf{e}-\overline{\mathbf{x}_2}\mathbf{x}_3)\mathbf{x}_3\big) ={\rm tr}(\mathbf{x}_2\overline{\mathbf{x}_3}){\rm tr}(\mathbf{x}_3^2\mathbf{x}_1)-{\rm tr}(\mathbf{x}_3^3\mathbf{x}_1\overline{\mathbf{x}_2})  \\
&={\rm tr}(\mathbf{x}_2\overline{\mathbf{x}_3}){\rm tr}((t_3\mathbf{x}_3\mathbf{x}_1-\mathbf{x}_1))
-{\rm tr}\big((t_3^2-1)\mathbf{x}_3\mathbf{x}_1\overline{\mathbf{x}_2}-t_3\mathbf{x}_1\overline{\mathbf{x}_2}\big) \\
&=(t_3t_{13}-t_1){\rm tr}(\mathbf{x}_2\overline{\mathbf{x}_3})-(t_3^2-1){\rm tr}(\mathbf{x}_3\mathbf{x}_1\overline{\mathbf{x}_2})
+t_3{\rm tr}(\mathbf{x}_1\overline{\mathbf{x}_2})   \\
&=(t_3t_{13}-t_1)(t_2t_3-t_{23})-(t_3^2-1)(t_2t_{13}-t_{123})+t_3(t_1t_2-t_{12}) \\
&=(t_3^2-1)t_{123}+t_1t_{23}+t_2t_{13}-t_3t_{12}-t_3t_{13}t_{23}.
\end{align*}
In the last line but three, we use $\mathbf{x}_3^2=t_3\mathbf{x}_3-\mathbf{e}$ and $\mathbf{x}_3^3=(t_3^2-1)\mathbf{x}_3-t_3\mathbf{e}$, which can be deduced by repeatedly applying (\ref{eq:H-C});
in the next line to last, we use ${\rm tr}(\mathbf{x}_i\overline{\mathbf{x}_j})=t_it_j-t_{ij}$.
Thus,
\begin{align*}
\xi&={\rm tr}(\mathbf{x}_2\mathbf{x}_3\overline{\mathbf{x}_1})+{\rm tr}(\mathbf{x}_3\overline{\mathbf{x}_2}\overline{\mathbf{x}_1})
+{\rm tr}(\mathbf{x}_3\mathbf{x}_1\overline{\mathbf{x}_2})+{\rm tr}(\mathbf{x}_3\mathbf{x}_1\overline{\mathbf{x}_3}\mathbf{x}_2\mathbf{x}_3) \\
&=(t_3^2-4)t_{123}-t_3t_{13}t_{23}+2t_1t_{23}+2t_2t_{13}.
\end{align*}

When $[\rho]\in\mathcal{X}_{3,i}$, we have $\xi=-4t_{123}$. Indeed,
\begin{itemize}
  \item[\rm(i)] if $[\rho]\in\mathcal{X}_{3,1}$, then $t_2=t_3=t_{23}=0$, so $\xi=-4t_{123}$;
  \item[\rm(ii)] if $[\rho]\in\mathcal{X}_{3,2}$, then $t_1=t_3=t_{13}=0$, so $\xi=-4t_{123}$;
  \item[\rm(iii)] if $[\rho]\in\mathcal{X}_{3,3}$, then $t_1=t_2=0$ and $t_3t_{123}=t_{13}t_{23}$, so $\xi=-4t_{123}$.
\end{itemize}
When $[\rho]\in\mathcal{X}_4$, from (\ref{eq:trace-expression}) we see
$$\xi=(t_3^2-4)\Big(\frac{\theta+1}{2}\Big)^2t_1t_2t_3-\Big(\frac{\theta+1}{2}\Big)^2t_1t_2t_3^3+4\Big(\frac{\theta+1}{2}\Big)t_1t_2t_3
=(1-\theta^2)t_1t_2t_3.$$

Therefore, (\ref{eq:TAP-2}) becomes
$$\Delta^{[\rho]}_B\doteq \prod_{j=1}^3(s_j+s_j^{-1}-t_j)
+\begin{cases}
0, & [\rho]\in\mathcal{X}^{\pm}_{1,i}\ \text{or\ }[\rho]\in\mathcal{X}_{2,i},  \\
4t_{123}+t_1t_2t_3, &[\rho]\in\mathcal{X}_{3,i}, \\
\theta^2t_1t_2t_3, &[\rho]\in\mathcal{X}_4. \end{cases}$$

\vspace{2mm}

\noindent
Haimiao Chen (orcid: 0000-0001-8194-1264), \emph{chenhm@math.pku.edu.cn}, \\
Tiantian Yu, \emph{she\_rasd@163.com}, \\
Depart of Mathematics, Beijing Technology and Business University, \\
Liangxiang Higher Education Park, Fangshan District, Beijing, China.

\end{document}